\documentclass[11pt]{article}

\usepackage{amsmath,amssymb}
\usepackage{latexsym}
\usepackage{graphicx}
\usepackage{bm}

\newtheorem{thm}{Theorem}

\newtheorem{rem}{Remark}

\newenvironment{proof}{\begin{trivlist}
                       \item[]{\bf Proof.}
                       \hspace{0cm}}{\hfill $\Box$
                       \end{trivlist}}

\def\diag{\mathop{\rm diag}}
\def\Re{\mathop{\rm Re}}

\begin{document}
\title{Dynamical systems gradient method for solving\\
ill-conditioned linear algebraic systems}

\author{N. S. Hoang$\dag$\footnotemark[1] \quad 
	  A. G. Ramm$\dag$\footnotemark[3] \\
\\
$\dag$Mathematics Department, Kansas State University,\\
Manhattan, KS 66506-2602, USA
}

\renewcommand{\thefootnote}{\fnsymbol{footnote}}
\footnotetext[1]{Email: nguyenhs@math.ksu.edu}
\footnotetext[3]{Corresponding author. Email: ramm@math.ksu.edu}

\date{}
\maketitle

\begin{abstract} \noindent A version of the Dynamical Systems Method 
(DSM) for
solving ill-conditioned linear algebraic systems is studied in this paper.  
An {\it a priori} and {\it a posteriori} stopping rules are justified.
An algorithm for computing the solution using a spectral decomposition of the
left-hand side matrix is proposed.  
Numerical results show
that when a spectral decompositon of the left-hand side matrix is available or not
computationally expensive to obtain the new method can be considered as an
alternative to the Variational Regularization.


\end{abstract}

\section{Introduction}

The Dynamical Systems Method (DSM) was systematically introduced and
investigated in \cite{R499} as a general method for solving operator
equations, linear and nonlinear, especially ill-posed operator equations.
In several recent publications various versions of the DSM, proposed in
\cite{R499}, were shown to be as efficient and economical as variational
regularization methods. This was demonstrated, for example, for the
problems of solving ill-conditioned linear algebraic systems \cite{R526},
 and stable numerical differentiation of noisy data
\cite{R415}, \cite{R441}, \cite{R530}.

The aim of this paper is to formulate a version of the DSM gradient method
for solving ill-posed linear equations and to demonstrate numerical
efficiency of this method. There 
is a large literature on
iterative regularization methods. These methods can be derived from 
a suitable version of the DSM by a discretization (see \cite{R499}). In 
the
Gauss-Newton-type version of the DSM one has to invert some linear
operator, which is an expensive procedure. The same is true for
regularized Newton-type versions of the DSM and of their iterative
counterparts. In contrast, the DSM gradient method we study in this 
paper {\it does not}
require inversion of operators.

We want to solve equation \begin{equation} \label{eqi1} Au=f,
\end{equation} where A is a linear bounded operator in a  Hilbert
space $H$. We assume that \eqref{eqi1} has a solution, possibly nonunique,
and denote by $y$ the unique minimal-norm solution to \eqref{eqi1},
$y\perp \mathcal{N}:=\mathcal{N}(A):= \{u: Au=0\}$, $Ay=f$. We assume that
the range of $A$, $R(A)$, is not closed, so problem \eqref{eqi1} is
ill-posed. Let $f_\delta$, $\|f-f_\delta\|\leq\delta$, be the noisy data.
We want to construct a stable approximation of $y$, given $\{\delta,
f_\delta, A\}$. There are many methods for doing this, see, e.g.,
\cite{I}, \cite{L}, \cite{M}, \cite{R499}, \cite{VV}, to mention a few
books, where variational regularization, quasisolutions, quasiinversion,
iterative regularization, and the DSM are studied.

The DSM version we study in this paper consists of solving the Cauchy problem
\begin{equation}
\label{eqi2}
\dot{u}(t)=-A^*(Au(t)-f),\quad u(0)=u_0, 
\quad u_0\perp N,\quad \dot{u}:=\frac{du}{dt},
\end{equation}
where $A^*$ is the adjoint to operator $A$, and proving the 
existence of the limit
$\lim_{t\to\infty}u(t)=u(\infty)$, and the relation $u(\infty)=y$, i.e.,
\begin{equation}
\label{eqi3}
\lim_{t\to\infty}\|u(t)-y\|=0.
\end{equation}

If the noisy data $f_\delta$ are given, then we solve the problem
\begin{equation}
\label{eqi4}
\dot{u}_\delta(t)=-A^*(Au_\delta(t)-f_\delta),\quad u_\delta(0)=u_0, 
\end{equation}
and prove that, for a suitable stopping time 
$t_\delta$, and $u_\delta :=u_\delta(t_\delta)$,
one has
\begin{equation}
\label{eqi5}
\lim_{\delta\to 0}\|u_\delta-y\|=0.
\end{equation}

In Section 2 these results are formulated precisely and 
recipes for choosing $t_\delta$ are proposed.

The novel results in this paper include the
proof of the discrepancy principle (Theorem 3), an efficient method
for computing $u_\delta(t_\delta)$ (Section \ref{sec comput}),
and an a priori stopping rule (Theorem 2). 

Our presentation is essentially self-contained.

\section{Results}

Suppose $ A: H\to H$ is a linear bounded operator in a Hilbert space $H$.
Assume that equation 
\begin{equation}
\label{eq1}
Au=f
\end{equation} 
has a solution not necessarily unique. Denote by $y$ 
the unique minimal-norm solution i.e., $y\perp \mathcal{N}:=\mathcal{N}(A)$. 
Consider the following Dynamical Systems Method (DSM)
\begin{equation}
\label{eq2}
\begin{split}
\dot{u}&=-A^*(Au-f),\\
u(0)&=u_0,
\end{split}
\end{equation}
where $u_0\perp \mathcal{N}$ is arbitrary. Denote $T:=A^*A,\, Q:= AA^*$.
The unique solution to \eqref{eq2} is
$$
u(t)=e^{-tT}u_0 + e^{-tT}\int_0^t e^{s T} ds A^*f.
$$
Let us show that any ill-posed linear equation \eqref{eq1} with exact data can be solved by the DSM.

\subsection{Exact data}

\begin{thm}
\label{thm1}
Suppose $u_0\perp \mathcal{N}$. Then problem \eqref{eq2} has a unique solution defined on
$[0,\infty)$, and $u(\infty)=y$, where $u(\infty)=\lim_{t\to\infty} u(t)$.
\end{thm}

\begin{proof}
Denote $w:=u(t)-y,\, w_0=w(0)$. Note that $w_0\perp \mathcal{N}$. One has
\begin{equation}
\label{eq3}
\dot{w}=-Tw,\quad T=A^*A.
\end{equation}
The unique solution to \eqref{eq3} is
$w=e^{-tT}w_0$. Thus,
$$
\|w\|^2=\int_0^{\|T\|} e^{-2t\lambda}d\langle E_\lambda w_0,w_0\rangle.
$$ 
where $\langle u,v\rangle$ is the inner product in $H$, and $E_\lambda$ is the resolution of the
identity of the selfadjoint operator $T$. Thus,
$$
\|w(\infty)\|^2 =\lim_{t\to \infty}\int_0^{\|T\|} e^{-2t\lambda}d\langle E_\lambda w_0,w_0\rangle=\|P_\mathcal{N}w_0\|^2=0,
$$
where $P_\mathcal{N}=E_0-E_{-0}$ is the orthogonal projector onto $\mathcal{N}$.
Theorem \ref{thm1} is proved.
\end{proof}

\subsection{Noisy data $f_\delta$}
\label{Sec2.2}

Let us solve stably equation \eqref{eq1} assuming that $f$ is not known,
but $f_\delta$, the noisy data, are known, 
where $\|f_\delta-f\|\le \delta$. Consider the following 
DSM 
$$
\dot{u}_\delta = - A^*(Au_\delta -f_\delta),\quad u_\delta(0)=u_0.
$$
Denote 
$$
w_\delta:=u_\delta-y,\quad T:=A^*A,\quad w_\delta(0)=w_0:=u_0-y\in 
\mathcal{N}^\perp.
$$ 
Let us prove the following result:

\begin{thm}
\label{thm2}
If $\lim_{\delta\to0}t_\delta = \infty,\,\lim_{\delta\to 0}t_\delta \delta=0$, and
$w_0\perp \mathcal{N}$, 
then 
$$
\lim_{\delta\to 0}\|w_\delta(t_\delta)\|=0.
$$
\end{thm}

\begin{proof}
One has
\begin{equation}
\label{eq4}
\dot{w}_\delta= -Tw_\delta + \eta_\delta,\quad \, \eta_\delta
=A^*(f_\delta-f),\quad \|\eta_\delta\|\le\|A\|\delta.
\end{equation}
The unique solution of equation \eqref{eq4} is
$$
w_\delta(t)=e^{-tT}w_\delta(0)+\int_0^te^{-(t-s)T}\eta_\delta ds.
$$
Let us show that 
$\lim_{t\to\infty} \|w_\delta(t)\|=0$. 
One has
\begin{equation}
\label{extra1}
\lim_{t\to\infty} \|w_\delta(t)\| \le \lim_{t\to\infty}\|e^{-tT}w_\delta(0)\|
+\lim_{t\to\infty}\bigg{\|}\int_0^te^{-(t-s)T}\eta_\delta ds\bigg{\|}.
\end{equation}
One uses the spectral theorem and gets:
\begin{equation}
\label{eq5}
\begin{split}
\int_0^te^{-(t-s)T}ds\eta_\delta&=\int_0^t\int_0^{\|T\|} dE_\lambda \eta_\delta e^{-(t-s)\lambda} ds\\
&=\int_0^{\|T\|} e^{-t\lambda}\frac{e^{t\lambda}-1}{\lambda}dE_\lambda\eta_\delta
=\int_0^{\|T\|}\frac{1-e^{-t\lambda}}{\lambda}dE_\lambda\eta_\delta.
\end{split}
\end{equation}
Note that
\begin{equation}
\label{eq6}
0\le\frac{1-e^{-t\lambda}}{\lambda}\le t,\quad \forall \lambda>0, t\ge 0,
\end{equation}
since $1-x\le e^{-x}$ for $x\ge 0$. 
From \eqref{eq5} and \eqref{eq6}, one obtains
\begin{equation}
\label{extra2}
\begin{split}
\bigg{\|}\int_0^te^{-(t-s)T}ds\eta_\delta\bigg{\|}^2
&=\int_0^{\|T\|}\big{|}\frac{1-e^{-t\lambda}}{\lambda}\big{|}^2d\langle E_\lambda\eta_\delta,\eta_\delta\rangle\\
&\le t^2\int^{\|T\|} d\langle E_\lambda\eta_\delta,\eta_\delta\rangle\\
&=t^2\|\eta_\delta\|^2.
\end{split}
\end{equation}
Since $\|\eta_\delta\|\le \|A\|\delta$, from \eqref{extra1} and \eqref{extra2}, one gets
$$
\lim_{\delta\to0} \|w_\delta(t_\delta)\| \le \lim_{\delta\to 0}\bigg{(}
\| e^{-t_\delta T}w_\delta(0)\|+t_{\delta}\delta\|A\|\bigg{)}=0.
$$
Here we have used the relation:
$$
\lim_{\delta\to 0}\| e^{-t_\delta T}w_\delta(0)\|=\|P_\mathcal{N}w_0\|=0,
$$
and the last equality holds because $w_0\in \mathcal{N}^\perp$. 
Theorem \ref{thm2} is proved.
\end{proof}

From Theorem \ref{thm2}, it follows that the relation $t_\delta=\frac{C}{\delta^\gamma}$, 
$\gamma=\text{const},\, \gamma\in(0,1)$ and $C>0$ is a constant, 
can be used as an \textit{a priori} stopping rule, i.e., for such $t_\delta$ one has
\begin{equation}\
\label{eq7}
\lim_{\delta\to0}\|u_\delta(t_\delta)-y\|=0.
\end{equation}

\subsection{Discrepancy principle}

Let us consider equation \eqref{eq1} with noisy data $f_\delta$, and a DSM of the form
\begin{equation}
\label{eq8}
\dot{u}_\delta = - A^*A u_\delta + A^*f_\delta,\quad u_\delta(0)=u_0.
\end{equation}
for solving this equation. 
Equation \eqref{eq8} has been used in Section~\ref{Sec2.2}. 
Recall that $y$ denotes the minimal-norm solution of equation \eqref{eq1}.

\begin{thm}
\label{thm3}
Assume that $\|Au_0-f_\delta\|> C\delta$. The
solution $t_\delta$ to the equation
\begin{equation}
\label{eq9}
h(t):=\|Au_\delta(t)- f_\delta\|=C\delta,\quad C=\text{const},\, C\in (1,2),
\end{equation}
does exist, is unique, and
\begin{equation}
\label{eq10}
\lim_{\delta\to 0} \|u_\delta(t_\delta)-y\|=0.
\end{equation}
\end{thm}

\begin{proof}
Denote 
$$
v_\delta(t):=Au_\delta(t)- f_\delta,\quad T:=A^*A,\quad Q=AA^*,\quad 
w(t):=u(t)-y,\quad w_0:=u_0-y.
$$ 
One has
\begin{equation}
\label{eq11}
\begin{split}
\frac{d}{dt}\|v_\delta(t)\|^2
&= 2\Re \langle A\dot{u}_\delta(t),Au_\delta(t)-f_\delta \rangle\\
&= 2 \Re \langle A[-A^*(Au_\delta(t) - f_\delta)],Au_\delta(t)-f_\delta 
\rangle\\
&=-2\|A^*v_\delta(t)\|^2\le 0.
\end{split}
\end{equation}
Thus, $\|v_\delta(t)\|$ is a nonincreasing function.
Let us prove that equation \eqref{eq9} has a solution for $C\in (1,2)$. 
Recall the known commutation formulas: 
$$
e^{-sT}A^*=A^*e^{-sQ},\, Ae^{-sT}=e^{-tQ}A.
$$
Using these formulas and the representation 
$$
u_\delta(t)=e^{-tT}u_0+\int_0^te^{-(t-s)T}A^*f_\delta ds,
$$ 
one gets:
\begin{align*}
v_\delta(t)
&=Au_\delta(t)-f_\delta\\
&=Ae^{-tT}u_0+A\int_0^te^{-(t-s)T}A^*f_\delta ds -f_\delta \\
&=e^{-t Q}Au_0+e^{-t Q}\int_0^{t} e^{sQ}dsQf_\delta-f_\delta \\
&=e^{-tQ}A(u_0-y)+e^{-tQ}f+e^{-tQ}(e^{tQ}-I)f_\delta - f_\delta\\
&=e^{-t Q}Aw_0 -e^{-t Q}f_\delta +e^{-t Q}f.
\end{align*}
Note that 
$$
\lim_{t\to\infty}e^{-t Q}Aw_0=\lim_{t\to\infty}Ae^{-t T}w_0 = AP_\mathcal{N}w_0=0.
$$ 
Here the continuity of $A$, and the following
relation
$$
\lim_{t\to\infty}e^{-tT}w_0=\lim_{t\to\infty}\int_0^{\|T\|}e^{-st}dE_sw_0=(E_0-E_{-0})w_0=P_\mathcal{N}w_0,
$$
were used. 
Therefore,
\begin{equation}
\label{eq12}
\lim_{t\to\infty}\|v_\delta(t)\|=\lim_{t\to\infty}\|e^{-t Q}(f-f_\delta)\|\le
\|f-f_\delta\|\le\delta,
\end{equation}
because $\|e^{-tQ}\|\le 1$. The function $h(t)$ is continuous on $[0,\infty)$, 
$h(0)=\|Au_0-f_\delta\|>C\delta$, $h(\infty)\le \delta$.
Thus, equation \eqref{eq9} must have a solution $t_\delta$.

Let us prove the uniqueness of $t_\delta$. Without loss of generality we can assume that there exists 
$t_1>t_\delta$ such that $\|Au_\delta(t_1)- f_\delta\|=C\delta$. Since $\|v_\delta(t)\|$ is 
nonincreasing and $\|v_\delta(t_\delta)\|=\|v_\delta(t_1)\|$, one has
$$
\|v_\delta(t)\|=\|v_\delta(t_\delta)\|,\quad \forall t\in [t_\delta, t_1].
$$
Thus,
\begin{equation}
\label{eq13}
\frac{d}{dt}\|v_\delta(t)\|^2=0,\quad \forall t\in (t_\delta, t_1).
\end{equation}
Using \eqref{eq11} and \eqref{eq13} one obtains
$$
A^*v_\delta(t)=A^*(Au_\delta(t)-f_\delta) = 0,\quad \forall t\in [t_\delta,t_1].
$$
This and \eqref{eq8} imply
\begin{equation}
\label{eq14}
\dot{u_\delta}(t)=0,\quad \forall t\in(t_\delta,t_1).
\end{equation}
One has
\begin{equation}
\label{eq15}
\begin{split}
\dot{u_\delta}(t)&=-Tu_\delta(t)+A^*f_\delta \\
&=-T\bigg{(}e^{-tT}u_0 +\int_0^te^{-(t-s)T}A^*f_\delta ds\bigg{)} + A^*f_\delta\\
&=-Te^{-tT}u_0 - (I-e^{-tT})A^*f_\delta+A^*f_\delta\\
&=-e^{-tT}(Tu_0 - A^*f_\delta).
\end{split}
\end{equation}
From \eqref{eq15} and \eqref{eq14}, one gets $Tu_0 - A^*f=e^{tT}e^{-tT}(Tu_0 - A^*f)=0$. 
Note that the operator $e^{tT}$ is an isomorphism for any fixed $t$ since $T$ is selfadjoint and bounded.
Since $Tu_0 - A^*f=0$, by \eqref{eq15} one has $\dot{u_\delta}(t)=0$,\, $u_\delta(t)=u_\delta(0)$,\, $\forall t\ge0$.
Consequently, 
$$
C\delta <\|Au_\delta(0)-f_\delta\|=\|Au_\delta(t_\delta)-f_\delta\|
=C\delta.
$$ 
This is a contradiction which proves
the uniqueness of $t_\delta$.

Let us prove \eqref{eq10}. First, we have the following estimate:
\begin{equation}
\label{eq16}
\begin{split}
\|Au(t_\delta)-f\|&\le \|Au(t_\delta)-Au_\delta(t_\delta)\|+\|Au_\delta(t_\delta)-f_\delta\|+\|f_\delta -f\|\\
&\le \bigg{\|}e^{-t_\delta Q}\int_0^{t_\delta}e^{sQ}Qds \bigg{\|} \|f_\delta-f\|+C\delta+\delta.
\end{split}
\end{equation}
Let us use the inequality:
$$
\big{\|}e^{-t_\delta Q}\int_0^{t_\delta}e^{sQ}Qds \big{\|}=\|I-e^{-t_\delta Q}\|\le 2,
$$
and conclude from \eqref{eq16}, that
\begin{equation}
\label{eq17}
\lim_{\delta\to0}\|Au(t_\delta)-f\|=0.
\end{equation}
Secondly, 
we claim that 
$$
\lim_{\delta\to0}t_\delta=\infty.
$$ 
Assume the contrary. Then there exist $t_0>0$ and a sequence
$(t_{\delta_n})_{n=1}^\infty$,
$t_{\delta_n}<t_0$, such that
\begin{equation}
\label{eq18}
\lim_{n\to\infty}\|Au(t_{\delta_n})-f\|=0.
\end{equation}
Analogously to \eqref{eq11}, one proves that 
$$
\frac{d\|v\|^2}{dt}\le 0,
$$ 
where $v(t):=Au(t)-f$.
Thus, $\|v(t)\|$ is nonincreasing. 
This and \eqref{eq18} imply the relation $\|v(t_0)\|=\|Au(t_0)-f\|=0$.
Thus,
$$
0=v(t_0)=e^{-t_0Q}A(u_0-y).
$$
This implies $A(u_0-y)=e^{t_0Q}e^{-t_0Q}A(u_0-y)=0$, so $u_0-y\in \mathcal{N}$. 
Since $u_0-y\in \mathcal{N}^\perp$, it follows that
 $u_0=y$. This is a contradiction because 
$$
C\delta\le\|Au_0-f_\delta\|=\|f-f_\delta\|\le\delta, \quad 1<C<2.
$$ 
Thus, $\lim_{\delta\to0}t_\delta=\infty$.

Let us continue the proof of \eqref{eq10}. Let 
$w_\delta(t):=u_\delta(t)-y$. 
We claim that $\|w_\delta(t)\|$ is nonincreasing on $[0,t_\delta]$. One has
\begin{align*}
\frac{d}{dt}\|w_\delta(t)\|^2&=2\Re \langle 
\dot{u_\delta}(t),u_\delta(t)-y\rangle\\
&= 2 \Re \langle -A^*(Au_\delta(t)-f_\delta),u_\delta(t)-y\rangle\\
&= -2 \Re \langle Au_\delta(t)-f_\delta,Au_\delta(t)-f_\delta 
+f_\delta-Ay\rangle\\
&\le -2 \|Au_\delta(t)-f_\delta\|\bigg{(}\|Au_\delta(t)-f_\delta\| - \|f_\delta -f\|\bigg{)}\\
&\le 0.
\end{align*}
Here we have used the inequalities: 
$$
\|Au_\delta(t)-f_\delta\|\ge C\delta >\|f_\delta -Ay\|=\delta,\quad  
\forall t\in [0,t_\delta].
$$ 

Let $\epsilon >0$ be arbitrary small. Since $\lim_{t\to\infty} u(t)=y$, there exists $t_0>0$, 
independent of $\delta$, such that
\begin{equation}
\label{eq19}
\| u(t_0)-y\|\le \frac{\epsilon}{2}.
\end{equation}
Since $\lim_{\delta\to0} t_\delta =\infty$, there exists $\delta_0$ such that 
$t_\delta>t_0,\, \forall \delta\in (0,\delta_0)$. 
Since $\|w_\delta(t)\|$ is nonincreasing 
on $[0,t_\delta]$ one has
\begin{equation}
\label{eqdel2}
\|w_\delta(t_\delta)\|\le \|w_\delta(t_0)\|\le \|u_\delta(t_0)-u(t_0)\|+ \|u(t_0)-y\|,\quad
\forall \delta\in(0,\delta_0).
\end{equation}
Note that
\begin{equation}
\label{eq21}
\begin{split}
\|u_\delta(t_0)-u(t_0)\|&=\|e^{-t_0T}\int_0^{t_0}e^{sT}ds A^*(f_\delta-f)\|
\le \|e^{-t_0T}\int_0^{t_0}e^{sT}dsA^*\|\delta.
\end{split}
\end{equation}
Since $e^{-t_0T}\int_0^{t_0}e^{sT}ds A^*$ is a bounded operator for any fixed $t_0$,
one concludes from \eqref{eq21} that $\lim_{\delta\to0}\|u_\delta(t_0)-u(t_0)\|=0$. 
Hence, there exists $\delta_1\in(0,\delta_0)$ such that
\begin{equation}
\label{eqdel3}
\|u_\delta(t_0)-u(t_0)\|\le \frac{\epsilon}{2},\quad \forall \delta\in(0,\delta_1).
\end{equation}
From \eqref{eq19}--\eqref{eqdel3}, one obtains
\begin{equation*}
\|u_\delta(t_\delta)-y\|=\|w_\delta(t_\delta)\|\le \frac{\epsilon}{2}+ \frac{\epsilon}{2}=\epsilon,\quad
\forall \delta\in(0,\delta_1).
\end{equation*}
This means that $\lim_{\delta\to0}u_\delta(t_\delta)=y$. Theorem~\ref{thm3} is proved.
\end{proof}

\section{Computing $u_\delta(t_\delta)$}
\label{sec comput}

\subsection{Systems with known spectral decomposition}

One way to solve the Cauchy problem \eqref{eq8} is 
to use explicit Euler or Runge-Kutta methods 
with a constant or adaptive stepsize $h$. 
However, stepsize $h$ for solving \eqref{eq8} by explicit numerical methods
is often smaller than 1 and the stopping time $t_\delta=nh$ may be large.
Therefore, the computation
time, characterized by the number of iterations $n$, for this approach may be large. 
This fact is also reported in
\cite{R526}, where one of the most efficient numerical methods for solving ordinary differential equations (ODEs),
the DOPRI45 (see \cite{Hairer}), is used for solving a Cauchy problem in a DSM. 
Indeed, the use of explicit Euler method leads to a Landweber iteration 
which is known for 
slow convergence.
Thus, it may be computationally expensive to compute $u_\delta(t_\delta)$ by
numerical methods for ODEs.

However, when $A$ in \eqref{eq8} is a matrix
and a decomposition $A=USV^*$, where $U$ and $V$ are unitary matrices
and $S$ is a diagonal matrix, is known, it is
possible to compute $u_\delta(t_\delta)$ at a speed comparable to 
other methods such as the variational regularization (VR) as it will be 
shown below.

We have
\begin{equation}
\label{eqforu}
u_\delta(t)=e^{-tT}u_0 + e^{-tT}\int_0^t e^{s T} ds A^*f_\delta,
\quad T:=A^*A.
\end{equation}
Suppose that a decomposition 
\begin{equation}
\label{specdecom}
A=USV^*,
\end{equation}
where $U$ and $V$ are unitary matrices
and $S$ is a diagonal matrix is known. 
These matrices possibly contain complex entries. Thus, $T=A^*A=V\bar{S}SV^*$ and $e^T=e^{V\bar{S}SV^*}$.
Using the formula $e^{V\bar{S}SV^*}=V e^{\bar{S}S}V^*$, which
is valid if $V$ is unitary and $\bar{S}S$ is diagonal,
equation~\eqref{eqforu} can be rewritten as
\begin{equation}
u_\delta(t)=Ve^{-t\bar{S}S}V^*u_0 + V\int_0^t e^{(s -t)\bar{S}S} ds \bar{S} U^*f_\delta.
\end{equation}
Here, the overbar stands for complex conjugation. 
Choose $u_0=0$. Then
\begin{equation}
\label{forforu}
u_\delta(t)=V\int_0^t e^{(s -t)\bar{S}S} ds \bar{S} h_\delta,\quad h_\delta:=U^*f_\delta.
\end{equation}

Let us assume that 
\begin{equation}
\label{eq32pre}
A^*f_\delta \not=0.
\end{equation} 
This is a natural assumption. Indeed, if $A^*f_\delta=0$, then by 
the definition of 
$h_\delta$ in \eqref{forforu}, relation $V^*V=I$, and equation \eqref{specdecom}, one gets 
\begin{equation}
\label{giaithich}
\bar{S}h_\delta = \bar{S}U^*f_\delta = V^*V\bar{S}U^*f_\delta = 
V^*A^*f_\delta=0.
\end{equation}
Equations \eqref{giaithich} and  \eqref{forforu} imply 
$u_\delta(t)\equiv 0$. 

The stopping time $t_\delta$ we choose by the following discrepancy 
principle:
$$
\|Au_{\delta}(t_\delta)-f_\delta\|=
\bigg{\|}\int_0^{t_\delta} e^{(s -t_\delta)\bar{S}S} ds \bar{S}S h_\delta-h_\delta\bigg{\|}
=\|e^{-t_\delta \bar{S}S}h_\delta\|=C\delta.
$$
where $1<C<2$.

Let us find $t_\delta$ from the equation
\begin{equation}
\label{eq32old}
\phi(t):=\psi(t)-C\delta =0,\qquad \psi(t):=\|e^{-t\bar{S}S}h_\delta\|.
\end{equation}
The existence and uniqueness of the solution $t_\delta$ to equation \eqref{eq32old} follow
from Theorem~\ref{thm3}. 

We claim that {\it equation \eqref{eq32old} can be solved by using Newton's iteration 
\eqref{eqnewton} for any initial value $t_0$ such that $\phi(t_0)>0.$ 
}

Let us prove this claim. It is sufficient to prove that
$\phi(t)$ is a monotone strictly convex function. This is proved below.

Without loss of generality, we can assume that $h_\delta$ (see \eqref{eq32old}) is a vector
with real components.
The proof remained essentially the same for $h_\delta$ with complex 
components.
%
%

First, we claim that 
\begin{equation}
\label{vocolen}
\sqrt{\bar{S}S}h_\delta\not=0,\qquad \text{and} \qquad
 \|\sqrt{\bar{S}S}e^{-t \bar{S}S}h_\delta\|\not=0,
\end{equation}
so $\psi(t)>0$.

Indeed, since $e^{-t \bar{S}S}$ is an isomorphism and $e^{-t \bar{S}S}$ commutes
with $\sqrt{\bar{S}S}$ one concludes that $\|\sqrt{\bar{S}S}e^{-t \bar{S}S}h_\delta\|=0$
iff $\sqrt{\bar{S}S}h_\delta=0$. If $\sqrt{\bar{S}S}h_\delta=0$ then $\bar{S}h_\delta=0$, and then, 
by equation ~\eqref{giaithich}, $A^*f_\delta=\bar{S}h_\delta=0$.
This contradicts to the assumption \eqref{eq32pre}. 

Let us now prove that $\phi$ monotonically decays and is strictly convex.
Then our claim will be proved.

One has
$$
\frac{d}{d t}\langle e^{-t \bar{S}S}h_\delta, e^{-t \bar{S}S}h_\delta\rangle
=-2\langle e^{-t\bar{S}S}h_\delta, \bar{S}S e^{-t \bar{S}S}h_\delta\rangle.
$$
Thus,
\begin{equation}
\label{eq32x}
\dot{\psi}(t)=\frac{d}{dt}\|e^{-t\bar{S}S}h_\delta\|
=\frac{\frac{d}{dt}\| e^{-t\bar{S}S}h_\delta \|^2}{2\| e^{-t\bar{S}S}h_\delta \|}
=-\frac{\langle e^{-t \bar{S}S}h_\delta, \bar{S}S e^{-t \bar{S}S}h_\delta\rangle}
{\|e^{-t \bar{S}S}h_\delta\|}.
\end{equation}
Equation \eqref{eq32x}, relation ~\eqref{vocolen}, and the fact that $\langle e^{-t \bar{S}S}h_\delta, \bar{S}S e^{-t \bar{S}S}h_\delta\rangle = \|\sqrt{\bar{S}S}e^{-t \bar{S}S}h_\delta\|^2$ imply
\begin{equation}
\label{eq36xc}
\dot{\psi}(t) < 0.
\end{equation}
%
From equation \eqref{eq32x} and the definition of $\psi$ in \eqref{eq32old}, one gets
\begin{equation}
\label{eq34x}
\psi(t)\dot{\psi}(t) = -\langle e^{-t \bar{S}S}h_\delta, \bar{S}S e^{-t \bar{S}S}h_\delta\rangle
\end{equation}
Differentiating equation \eqref{eq34x} with respect to $t$, one obtains
\begin{align*}
\psi(t)\ddot{\psi}(t) + \dot{\psi}^2(t) &= \langle \bar{S}Se^{-t \bar{S}S}h_\delta, \bar{S}S e^{-t \bar{S}S}h_\delta\rangle + 
\langle e^{-t \bar{S}S}h_\delta, \bar{S}S\bar{S}S e^{-t \bar{S}S}h_\delta\rangle\\
&= 2 \|\bar{S}Se^{-t \bar{S}S}h_\delta\|^2.
\end{align*}
This equation and equation \eqref{eq32x} imply
\begin{equation}
\label{eq38xc}
\psi(t)\ddot{\psi}(t) = 2 \|\bar{S}Se^{-t \bar{S}S}h_\delta\|^2 - 
\frac{\langle e^{-t \bar{S}S}h_\delta, \bar{S}S e^{-t \bar{S}S}h_\delta\rangle^2}
{\|e^{-t \bar{S}S}h_\delta\|^2} \ge \|\bar{S}Se^{-t \bar{S}S}h_\delta\|^2 > 0.
\end{equation}
Here the inequality: $\langle e^{-t \bar{S}S}h_\delta, 
\bar{S}S e^{-t \bar{S}S}h_\delta\rangle\le
\|e^{-t \bar{S}S}h_\delta\|\|\bar{S}S e^{-t \bar{S}S}h_\delta\|$ was used.
Since $\psi>0$, inequality \eqref{eq38xc} implies
\begin{equation}
\label{eq39xc}
\ddot{\psi}(t) > 0.
\end{equation}
It follows from inequalities \eqref{eq36xc} and \eqref{eq39xc} that $\phi(t)$ 
is a strictly convex and decreasing function on $(0,\infty)$. 
Therefore, $t_\delta$ can be found by Newton's iterations:
\begin{equation}
\label{eqnewton}
\begin{split}
t_{n+1}&= t_n -\frac{\phi(t_n)}{\dot{\phi}(t_n)}\\
&=t_n +
\frac{\| e^{-t_n
\bar{S}S}h_\delta\| - C\delta}{\langle \bar{S}Se^{-t_n
\bar{S}S}h_\delta,e^{-t_n \bar{S}S}h_\delta\rangle}
\| e^{-t_n\bar{S}S}h_\delta\|,\quad n=0,1,...,
\end{split}
\end{equation}
for any initial guess $t_0$ of $t_\delta$ such that $\phi(t_0)>0$. 
Once $t_\delta$ is found, the solution $u_\delta(t_\delta)$ is computed by \eqref{forforu}.

\begin{rem}{\rm
In the decomposition $A = VSU^*$ we do not assume that $U,V$ and $S$ are matrices
with real entries. The singular value decomposition (SVD) is a particular
case of this decomposition. 

It is computationally expensive to get the SVD of a matrix in general. 
However, there are many problems in which the decomposition 
\eqref{specdecom} can be computed
fast using the fast Fourier transform (FFT). Examples include image restoration
problems with circulant block matrices (see \cite{N}) 
and deconvolution problems.
(see Section~\ref{image sec}). 
}
\end{rem}

\subsection{On the choice of $t_0$}

Let us discuss a strategy for choosing the initial value $t_0$ in Newton's iterations
for finding $t_\delta$.
We choose $t_0$ satisfying condition:
\begin{equation}
\label{eqt0}
0<\phi(t_0)=\|e^{-t_0\bar{S}S}h_\delta\|-\delta \le\delta
\end{equation}
by the following strategy
\begin{enumerate}
\item{Choose $t_0:=10\frac{\|h_\delta\|}{\delta}$ as an initial guess for $t_0$.}
\item{Compute $\phi(t_0)$. If $t_0$ satisfying \eqref{eqt0} we are done. Otherwise, we go to step 3.}

\item{If $\phi(t_0)< 0$ and the inequality $\phi(t_0)> \delta$ has not occurred in iteration, 
we replace $t_0$ by $\frac{t_0}{10}$  and go back to step 2. 
If $\phi(t_0)< 0$ and the inequality $\phi(t_0)> \delta$ has occurred in iteration, 
we replace $t_0$ by $\frac{t_0}{3}$ and go back to step 2.
If $\phi(t_0)> \delta$, we go to step 4.}


\item{
If $\phi(t_0)> \delta$ and the inequality $\phi(t_0)< 0$ has not occured in iterations, 
 we replace $t_0$ by $3t_0$ and go back to step 2.
If the inequality $\phi(t_0)< 0$ has occured in some iteration before, we stop the iteration and use $t_0$ 
as an initial guess in  
Newton's iterations for finding $t_\delta$.}
\end{enumerate}


\section{Numerical experiments}

In this section results of some numerical experiments with ill-conditioned linear algebraic systems
are reported.
In all the experiments, by DSMG we denote the version of the DSM described in this paper, 
by VR we denote the Variational Regularization,
implemented using the discrepancy principle, and by DSM-\cite{R526} we denote 
the method developed in \cite{R526}.

\subsection{A linear algebraic system for the computation of second derivatives}

Let us do some numerical experiments with linear algebraic systems arising in a
numerical experiment of computing the second derivative of a noisy function.

The problem is reduced to an integral equation of the first kind.
A linear algebraic system is obtained by a discretization of the integral equation whose kernel $K$ is Green's function
$$
K(s,t)=
\left\{
\begin{matrix}
s(t-1),\quad \text{if}\quad s<t\\
t(s-1),\quad \text{if}\quad s\geq t
\end{matrix}
\right. .
$$
Here $s,t\in[0,1]$.  
Using $A_N$ from \cite{R526}, we do some numerical experiments for solving $u_N$ from the linear algebraic system
$A_N u_N=b_{N,\delta}$. In the experiments
the exact right-hand side is computed by the formula $b_N=A_Nu_N$ when $u_N$ is given. 
In this test, $u_N$ is computed by 
$$
u_N:=\big{(}u(t_{N,1}),u(t_{N,2}),....,u(t_{N,N})\big{)}^T,\qquad t_{N,i}:=\frac{i}{N},\quad i=1,...,N,
$$
where $u(t)$ is a given function. 
We use $N=10,20,...,100$ and $b_{N,\delta} = b_N + e_N$, 
where $e_N$ is a random vector whose coordinates are independent, normally distributed, with mean 0 and  variance 1, 
and scaled so that $\|e_N\|=\delta_{rel}\|b_N\|$. 
This linear algebraic system is mildly ill-posed: the condition 
number of $A_{100}$ is $1.2158\times10^{4}$. 

\begin{figure}[!h!t!b!]
\centerline{
\begin{tabular}{c}
\mbox{\includegraphics[scale=0.92]{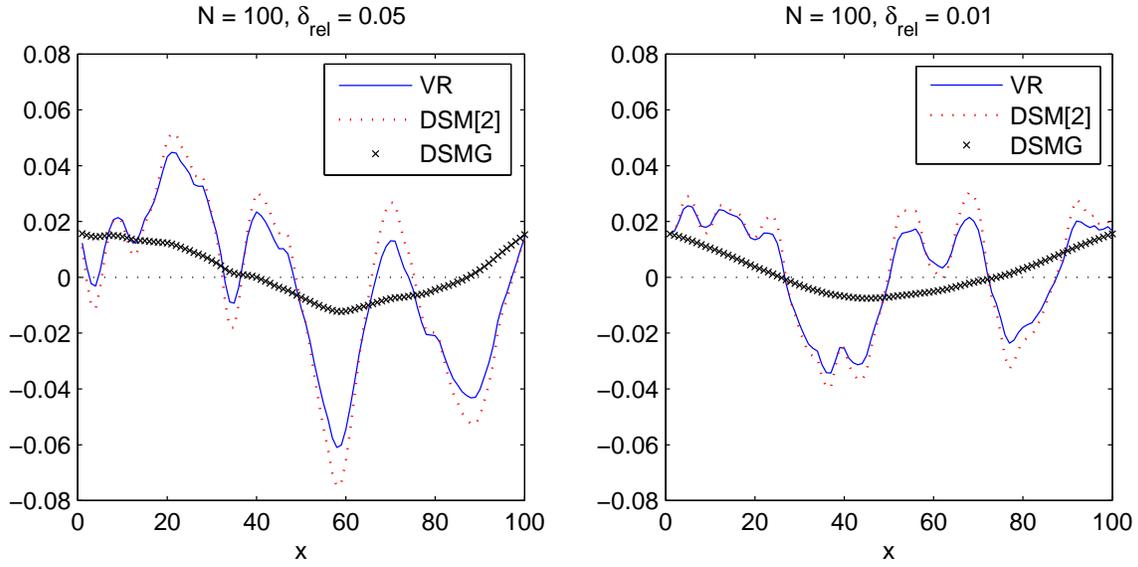}}
\end{tabular}
}
\caption{\it Plots of differences between the exact solution and solutions obtained
by the DSMG, VR and DSM-\cite{R526}.
}
\label{fig1}
\end{figure}

In Figure~\ref{fig1}, the difference between the exaction and
solution obtained by the DSMG, VR and DSM-\cite{R526} are plotted. 
In these experiments, we used $N=100$ and $u(t)=\sin(\pi t)$ with $\delta_{rel}=0.05$
and $\delta_{rel}=0.01$. 
Figure~\ref{fig1} shows that the results obtained by
the VR and the DSM-\cite{R526} are very close to each other.
The results obtained by the DSMG are much better than those by 
the DSM-\cite{R526} and by the VR. 

Table~\ref{table2} presents numerical results when $N$ varies from $10$ to $100$,
$u(t)=\sin(2\pi t)$, and $t\in[0,1]$.
In this experiment the DSMG yields more accurate solutions than the 
DSM-\cite{R526} and the VR. The DSMG in this experiment takes more iterations than the 
DSM-\cite{R526} and the VR to get a solution. 
\begin{table}[h] 
\caption{Numerical results for computing second derivatives with $\delta_{rel}=0.01$.}
\label{table2}
\centering
\small
\begin{tabular}{@{  }c@{\hspace{2mm}}|c@{\hspace{2mm}}
@{\hspace{2mm}}c@{\hspace{2mm}}|c@{\hspace{2mm}}c@{\hspace{2mm}}
|c@{\hspace{2mm}}c@{\hspace{2mm}}|l@{}} 
\hline
&
\multicolumn{2}{c|}{DSM}
&\multicolumn{2}{c|}{DSM-\cite{R526}}&\multicolumn{2}{c|}{VR}\\
$N$&
n$_{\text{iter}}$&$\frac{\|u_\delta-y\|_{2}}{\|y\|_2}$&
n$_{\text{linsol}}$&$\frac{\|u_\delta-y\|_{2}}{\|y\|_2}$&
n$_{\text{linsol}}$&$\frac{\|u_\delta-y\|_{2}}{\|y\|_2}$\\
\hline
20    &9    &0.0973    &3    &0.1130    &6    &0.1079 \\
30    &5    &0.0831    &4    &0.1316    &6    &0.1160 \\
40    &7    &0.0488    &4    &0.1150    &6    &0.1045 \\
50    &9    &0.0614    &4    &0.1415    &6    &0.1063 \\
60    &6    &0.0419    &4    &0.0919    &6    &0.0817 \\
70    &9    &0.0513    &4    &0.0961    &6    &0.0842 \\
80    &6    &0.0418    &4    &0.1225    &6    &0.0981 \\
90    &7    &0.0287    &4    &0.0919    &7    &0.0840 \\
100   &7    &0.0248    &5    &0.0778    &7    &0.0553 \\
\hline 
\end{tabular}
\end{table}

In this experiment the DSMG is implemented using the SVD of $A$ obtained by the function {\it svd} in Matlab.
As already mentioned, the SVD is a special case of the spectral decomposition 
\eqref{specdecom}. 
It is expensive to compute the SVD, in general. However, there are practically important problems where the spectral decomposition
\eqref{specdecom} can be computed fast (see Section~\ref{image sec} below). 
These problems consist of deconvolution problems using the Fast Fourier Transform (FFTs).

The conclusion from this experiment is: the DSMG may yield results with much better accuracy than the VR and DSM-\cite{R526}.
Numerical experiments for various $u(t)$ show that the DSMG competes favorably with the VR and the DSM-\cite{R526}.

\subsection{An application to image restoration}
\label{image sec}

The image degradation process can be modeled by the following equation:
\begin{equation}
\label{eq35}
g_\delta = g+w,\quad g = h\ast f,\quad \|w\|\le \delta,
\end{equation}
where $h$ represents a convolution function that models the blurring that many imaging systems introduce. 
For example, camera defocus, motion blur, imperfections of the lenses, all 
these phenomenon can be modeled by choosing a suitable $h$. 
The functions $g_\delta$, $f$, and $w$ are the observed image, the original signal, and the noise, respectively. 
The noise $w$ can be due to the electronics used (thermal and shot noise), the recording medium (film grain), 
or the imaging process (photon noise). 

In practice $g,h$ and $f$ in equation \eqref{eq35} are often given as functions
 of a discrete argument and equation
\eqref{eq35} can be written in this case as
\begin{equation}
g_{\delta,i} = g_i+w_i = \sum_{j=-\infty}^\infty f_j h_{i-j}+w_i,\quad i\in \mathbb{Z}.
\end{equation}
Note that one (or both) signals $f_j$ and $h_j$ have compact support (finite length). 
Suppose that signal $f$ is periodic with period $N$, i.e., $f_{i+N}=f_i$, and $h_j=0$
 for $j< 0$ and $j\ge N$.  
Assume that $f$ is represented by a sequence $f_0,...,f_{N-1}$
and $h$ is represented by $h_0,...,h_{N-1}$. Then the convolution $h\ast f$ is periodic signal $g$ with period $N$, and the elements of $g$ are defined as
\begin{equation}
\label{eq36}
g_i = \sum_{j=0}^{N-1} h_j f_{(i-j)\, mod\, N},\quad i=0,1,...,N-1.
\end{equation}
Here $(i-j)\, mod\, N$ is $i-j$ modulo $N$. 
The discrete Fourier transform (DFT) of $g$ is defined as the sequence
$$
\hat{g}_k:=\sum_{j=0}^{N-1}g_je^{-i2\pi jk/N},\qquad k=0,1,...,N-1.
$$
Denote $\hat{g}=(\hat{g}_0,....,\hat{g}_{N-1})^T$. Then equation \eqref{eq36} implies
\begin{equation}
\label{eq37}
\hat{g} = \hat{f}\hat{h},\qquad \hat{f}\hat{g}:=(\hat{f}_0\hat{h}_0,\hat{f}_1\hat{h}_1,...,\hat{f}_{N-1}\hat{h}_{N-1})^T.
\end{equation}
Let $\diag(a)$ denote a diagonal matrix whose diagonal is $(a_0,...,a_{N-1})$ and other entries are zeros. Then equation \eqref{eq37} can be rewritten as
\begin{equation}
\label{eq38}
\hat{g} = A\hat{f},\qquad A:=\diag(\hat{h}).
\end{equation}
Since $A$ is of the form \eqref{specdecom} with $U=V=I$
and $S=\diag(\hat{h})$,
 one can use the DSMG method to solve equation \eqref{eq38} stably for $\hat{f}$.

\begin{figure}[!h!t!b!]
\centerline{
\begin{tabular}{c}
\mbox{\includegraphics[scale=1]{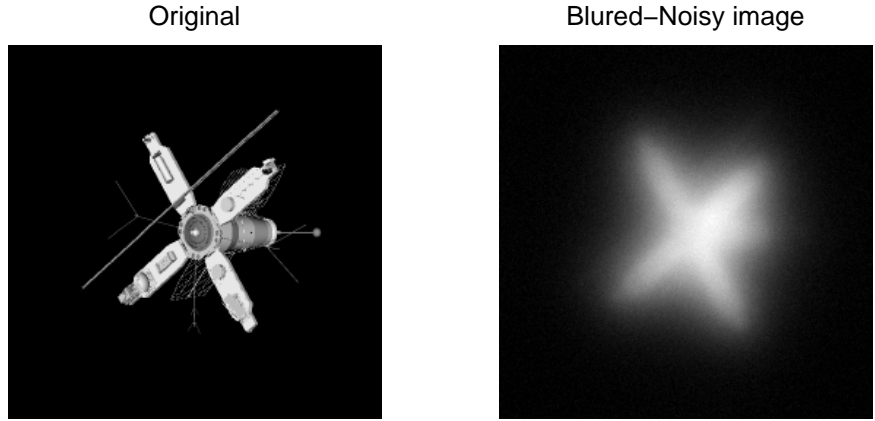}}
\end{tabular}
}
\caption{\it Original and Blurred-noisy images. 
}
\label{figimage1}
\end{figure}

The image restoration test problem we use is taken from \cite{N}.
This test problem was developed at the US Air
Force Phillips Laboratory, Lasers and Imaging Directorate, Kirtland Air Force
Base, New Mexico. 
The original and blurred
images have $256\times 256$ pixels, and are shown in Figure~\ref{figimage1}.
These data has been widely used in the literature for testing image restoration
algorithms.

\begin{figure}[!h!t!b!]
\centerline{
\begin{tabular}{c}
\mbox{\includegraphics[scale=1]{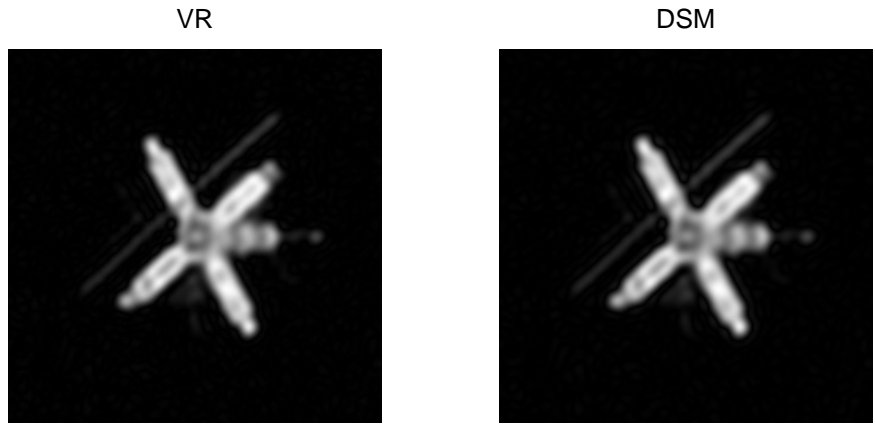}}
\end{tabular}
}
\caption{\it Regularized images when noise level is 1\%. 
}
\label{figimage2}
\end{figure}

Figure~\ref{figimage2} plots the regularized images by the VR and the DSMG when $\delta_{rel}=0.01$. Again, with an input value for $\delta_{rel}$, the observed blurred-noisy
images is computed by
$$
g_\delta = g + \delta_{rel}\frac{\|g\|}{\|err\|}err,
$$
where $err$ is a vector with random entries normally distributed with mean 0 and variance 1. 
In this experiment, it took 5 and 8 iterations for the DSMG and the VR, respectively, to yield numerical results.
From Figure~\ref{figimage2} one concludes that the DSMG is comparable to the VR in terms of
accuracy. The time of computation in this experiment is about the same for the VR and DSMG.

\begin{figure}[!h!t!b!]
\centerline{
\begin{tabular}{c}
\mbox{\includegraphics[scale=1]{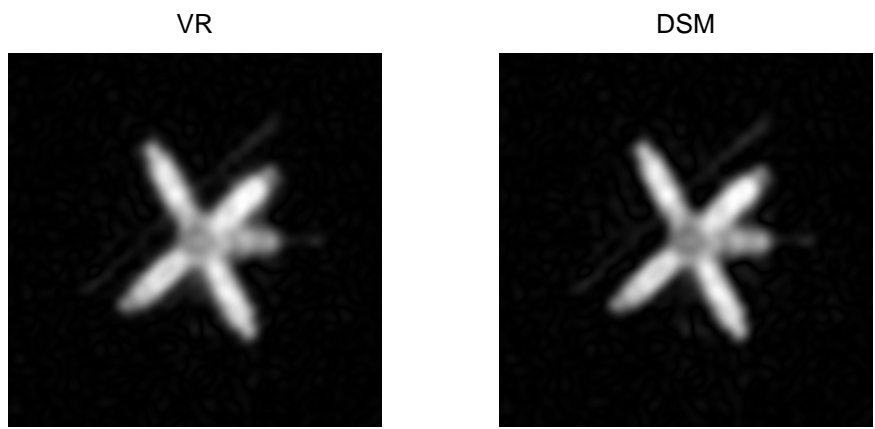}}
\end{tabular}
}
\caption{\it Regularized images when noise level is 5\%. 
}
\label{figimage3}
\end{figure}

Figure~\ref{figimage3} plots the regularized images by the VR and the DSMG when $\delta_{rel}=0.05$.
It took 4 and 7 iterations for the DSMG and the VR, respectively, to yield numerical results.
Figure~\ref{figimage3} shows that the images obtained by the DSMG and the VR are
about the same.

The conclusions from this experiment are: the DSMG yields results with
the same accuracy as the VR, and requires less iterations than the VR. 
The restored images by
the DSM-\cite{R526} are about the same as those by the VR.

\begin{rem}{\rm
Equation \eqref{eq35} can be reduced to equation \eqref{eq37} whenever one of
the two functions $f$ and $h$ has compact support and the other is periodic.
}
\end{rem}

\section{Concluding remarks}

A version of the Dynamical Systems Method for solving ill-conditioned linear algebraic systems is studied in this paper. 
An {\it a priori} and {\it a posteriori} stopping rules are formulated and justified. 
An algorithm for computing the solution
in the case when a spectral decomposition of the matrix $A$ is available is presented. 
Numerical results show that the DSMG, i.e., the DSM version developed in this paper, yields results comparable to those obtained by the 
VR and the DSM-\cite{R526} developed in \cite{R526}, and the DSMG method may yield much more accurate results than the VR method. 
It is demonstrated in \cite{N} that
the rate of convergence of the Landweber method can be increased by 
using preconditioning techniques.
The rate of convergence of the DSM version, presented in this paper, 
might be improved by a similar technique. 
The advantage of our method over the steepest descent in \cite{N} is the following:
{\it the stopping time $t_\delta$ can be found 
from a discrepancy principle by Newton's iterations for a wide range of  initial guess
$t_0$; 
when $t_\delta$ is found one can compute the solution without any iterations.}
Also, our method requires less iterations than the steepest descent in \cite{N}, which is an accelerated version of the Landweber method.

\end{document}